\documentclass[11pt]{article}
\newtheorem{theorem}{Theorem}
\newtheorem{lemma}[theorem]{Lemma}
\newtheorem{prop}[theorem]{Proposition}

\newcommand{\nl}{\newline}
\newcommand{\dist}{{\rm dist}}
\newcommand{\N}{{\bf N}}

\newcommand{\R}{{\bf R}}

\newcommand{\cG}{{\cal G}}
\newcommand{\cD}{{\cal D}}

\newcommand{\cA}{{\cal A}}

\newcommand{\D}{{\rm (D1)}}
\newcommand{\DD}{{\rm (D2)}}
\newcommand{\DDD}{{\rm (D3)}}

\newcommand{\po}{{\partial\Omega}}
\newcommand{\dom}{{\rm Dom}}

 \newcommand{\inprod}[2]{{\langle{#1},{#2}\rangle}}

 \newcommand{\darr}[4]{{\left\{\begin{array}{ll}
   {#1}&{#2}\\
   {#3}&{#4}
 \end{array}\right.}}

\newcommand{\ale}{{\rm a.e.}}
\newcommand{\ia}{({\rm i})}
\newcommand{\ib}{({\rm ii})}

\newcommand{\sscr}[2]{\scriptstyle {#1} \atop\scriptstyle {#2}}
\newcommand{\ssscr}[3]{\scriptstyle {#1} \atop{\scriptstyle {#2}\atop\scriptstyle {#3}}}

\newcommand{\pa}{{\bf (P_{\alpha})}}

\title{Boundary decay estimates for solutions of fourth-order elliptic equations} 

\author{
G. Barbatis\\  
Department of Mathematics \\
University of Ioannina\\ 
45110 Ioannina, Greece
}
\date{}

\begin{document}

\maketitle

\begin{abstract}
We obtain integral boundary decay estimates for solutions of fourth-order elliptic equations
on a bounded domain with regular boundary. We apply these estimates to obtain stability
bounds for the corresponding eigenvalues under small perturbations of the boundary.

\vspace{0.3cm}

\end{abstract}

\section{Introduction}\label{intro}

\parskip=0pt

\

Let $\Omega$ be a bounded region in $\R^N$ and let
$H$ be a fourth-order, self-adjoint, uniformly elliptic operator on $L^2(\Omega)$
subject to Dirichlet boundary conditions on $\po$,
\[
Hu(x)=\sum_{\sscr{\mid\eta\mid\leq 2}
{\mid\zeta\mid\leq 2}}
D^{\eta}\{a_{\eta\zeta}(x)D^{\zeta}u(x)\}\; , \quad x\in\Omega\, .
\]
The scope of this paper is to obtain integral boundary decay estimates for solutions of the equation
\begin{equation}
Hu=f\; , \quad f\in L^2(\Omega)\, .
\label{mera}
\end{equation}
More precisely, we want to establish ranges of $\beta>0$ for which the integrals
\[ \int d^{-2-\beta}u^2dx \quad \mbox{ and } \quad \int d^{-\beta}|\nabla u|^2dx
\]
($d(x)=\dist(x,\po)$) are finite.
If the boundary $\Omega$ is regular in the sense that the Hardy-Rellich inequality
\[
\int_{\Omega}(\Delta u)^2dx \geq c\int_{\Omega}\Big\{\frac{|\nabla u|^2}{d^2}+\frac{u^2}{d^4}\Big\}dx\; , \quad
u\in H^2_0(\Omega)\, ,
\]
is valid, we then immediately have such an estimate since $H^2_0(\Omega)=\dom(H^{1/2})$. Our aim is to establish
better decay estimates that exploit the fact that the solution $u$ of (\ref{mera}) belongs not only in $H^2_0(\Omega)$
but also in $\dom(H)$.

This problem is well studied in the case of second-order operators.
In \cite{Dmz} Davies obtained boundary decay estimates of the form
\begin{equation}
\int_{\Omega}\left( \frac{|\nabla u|^2}{d^{2\alpha}} +
\frac{u^2}{d^{2+2\alpha}}\right)dx \leq c\biggl(\|Hu\|_2\|H^{1/2}u\|_2 + \|u\|_2^2\biggr), 
\quad u\in\dom(H),
\label{ww1}
\end{equation}
for $\alpha>0$ in some interval $(0,\alpha_0)$.
Here $\alpha_0$ is an explicitly given constant
which depends on the boundary regularity and the ellipticity constants of $H$.
As an application of (\ref{ww1}) stability estimates were obtained on the eigenvalues
$\{\lambda_n\}$ of $H$ under small perturbations of the boundary $\po$.
More precisely, if $\tilde{\Omega}\subset\Omega$ is a domain such that
$\partial\tilde{\Omega}\subset\{x\in\Omega \, : \, \dist(x,\po)<\epsilon\}$
and if $\{\tilde{\lambda}_n\}$ are the corresponding Dirichlet eigenvalues (the operator $\tilde{H}$ being
defined by form restriction), then it was shown that (\ref{ww1}) implies
\begin{equation}
0\leq\tilde{\lambda}_n-\lambda_n\leq c_n\epsilon^{2\alpha} 
\label{in}
\end{equation}
for all $n\in\N$ and all $\epsilon>0$ small enough. This estimate has obvious applications in the numerical computation of eigenvalues;
see \cite{Dmz} for more on specific examples.

Inequality (\ref{in}) was subsequently improved in \cite{Dcam}, where $\epsilon^{2\alpha}$, $\alpha<\alpha_0$,
was replaced by $\epsilon^{2\alpha_0}$, for the same $\alpha_0$; this was done by estimating the integrals
$\int_{d(x)<\epsilon}|\nabla u|^2dx$ and $\int_{d(x)<\epsilon}u^2dx$ for small $\epsilon>0$. 
For results analogous to those of \cite{Dcam} for the $p$-Laplacian together with
applications we refer to Fleckinger et al. \cite{FHT}. See also \cite{EHK} where estimates of this type were
first obtained for eigenfunctions of second-order operators. For relevant results in the case of singular operators
see \cite{M}.

%%%%%%%%%%%%%%%%%%%%%

In our main theorem we obtain integral decay estimates analogous to (\ref{ww1}) for fourth-order operators.
More precisely, for a fourth-order operator
$H$ with $L^{\infty}$ coefficients we establish boundary decay estimates of the form
\begin{equation}
\int_{\Omega}\left( \frac{|\nabla^2u|^2}{d^{2\alpha}}+\frac{|\nabla u|^2}{d^{2+2\alpha}} +
\frac{u^2}{d^{4+2\alpha}}\right)dx \leq c\biggl(\|Hu\|_2\|H^{\alpha/2}u\|_2 + \|u\|_2^2\biggr), 
\;\; u\in\dom(H),
\label{elv}
\end{equation}
for $\alpha$ in an interval $(0,\alpha_0)$. Under additional assumptions we obtain $\alpha_0=1/2$,
which is optimal. To prove (\ref{elv}) we first use some general inequalities, which lead to a property
$\pa$ being identified as sufficient for the validity of (\ref{elv}). We then study property $\pa$, and find sufficient
conditions under which it is valid; the distance function used here is taken to be the Finsler distance
induced by the operator.

Technical reasons oblige us to make a regularity assumption that is not needed in the second-order case and requires
$d(x)$ to be $C^2$ near $\po$. This relates to a recurrent and largely unsolved issue in higher-order
problems: a distance function
is normally only once differentiable, but is required for technical reasons to be differentiated more than once;
see for example \cite{B}, where such an issue has arisen in the context of heat kernel estimates.

Finally, as an application of (\ref{elv}) we obtain stability bounds analogous to (\ref{in})
on the eigenvalues of $H$ under small boundary perturbations.

The structure of the paper is as follows: in Section \ref{ageneral} we provide a sufficient
condition $\pa$ for the validity of (\ref{elv}); in Section \ref{theproperty}
we establish the range of $\alpha$ for which $\pa$ is valid for different classes of operators;
and in Section \ref{anapplication} we present the application to eigenvalue stability.

\parindent=0pt
\parskip=5pt

\subsection*{Setting}
We fix some notation. Given a multi-index $\eta=(\eta_1,\ldots,\eta_N)$ we write
$\eta!=\eta_1!\ldots\eta_N!$ and $|\eta|=\eta_1+\cdots+\eta_N$.
We write $\gamma\leq\eta$ if $\gamma_i\leq \eta_i$ for all $i$, in which case we
also set $c^{\eta}_{\gamma}=\eta !/\gamma ! (\eta -\gamma)!$.
We use the standard notation
$D^{\eta}u=(\partial/\partial x_1)^{\eta_1}\ldots (\partial/\partial x_N)^{\eta_N}u$
and $(\nabla u)^{\eta}=u_{x_1}^{\eta_1}\ldots u_{x_N}^{\eta_N}$.
By $\nabla^2u$ we denote the vector $(u_{x_ix_j})_{i,j=1}^N$. The letter $c$ will denote a
constant whose value may change from line to line; the constants $c_1,c_2$ and $c_3$ however
are the same throughout the paper.

We now describe our setting. We assume that $\Omega$ is a bounded domain in $\R^N$ with
boundary $\po$. We consider a distance function $d(\cdot ,\cdot)$ on $\Omega$,
and denote by $d(\cdot)$ the corresponding distance to the boundary $\po$.
We say that $d(\cdot)$ belongs in the class $\cD$ if it satisfies:
\begin{eqnarray*}
\D && \mbox{There exist $c_1,c_2>0$ such that for any $x,y\in\Omega$}\\
&& \quad {\rm (a)} \qquad  c_1 \leq |\nabla_z d(z,y)| \leq c_2 \, , \quad z\in\Omega \\
&& \quad {\rm (b)} \qquad  c_1 d_{Euc}(x,y)\leq d(x,y)\leq c_2 d_{Euc}(x,y) .\\[0.2cm]
 \DD && \mbox{There exist $\theta,\tau >0$ such that} \\
&&\quad {\rm (a)} \qquad \mbox{$d(x)$ is $C^2$ on $\{0<d(x)<\theta\}$} \\
&&\quad {\rm (b)} \qquad \mbox{$|\nabla^2d(x)|\leq cd(x)^{-1+\tau}$ on $\{0<d(x)<\theta\}$.} \\[0.2cm]
\DDD &&\mbox{The following Hardy-Rellich inequalities are valid for all $v\in C^2_c(\Omega)$:}\\
&&\quad {\rm (a)} \qquad \int_{\Omega}(\Delta v)^2dx \geq c_3
\int_{\Omega} \frac{v^2}{d^4}dx  \\
&&\quad {\rm (b)} \qquad \int_{\Omega}(\Delta v)^2dx \geq c_3
\int_{\Omega} \frac{|\nabla v|^2}{d^2}dx \, .
\end{eqnarray*}
We note that a sufficient condition for $\DDD$(b) is the Hardy inequality
\[
\int_{\Omega}|\nabla v|^2dx \geq c_3\int_{\Omega}\frac{v^2}{d^2}dx \, .
\]
The distance $d(\cdot,\cdot)$ will typically be a Finsler distance, in which case (a) and (b) of $\D$ are
equivalent.
Condition $\DD$ is essentially a strong regularity assumption on $\po$, as will be seen below.
Its validity in
examples will always involve the specific value $\tau=1$; we choose however this more general and somewhat axiomatic
setting because, we believe, it shows more clearly what the essential ingredients are.
Finally, for more on Hardy-Rellich inequalities, optimal constants as well as improved versions of
such inequalities we refer to \cite{BFT,BT} and references therein.

In the sequel we shall often need to twice differentiate $d(x)$ near $\po$.
In order to avoid repeatedly
splitting integrals in two, we redefine $d(x)$ on $\{x\in\Omega : d(x)>\theta\}$ so that now
$d(x)$ is $C^2$ and positive function on $\Omega$ such that
$d(x)$ equals $\inf\{d(x,y) : y\in\po\}$
for $x\in\{\dist(x,\po)<\theta\}$ but not necessarily for all $x\in\Omega$ (of course $d(x,\cdot)$
extends to $\po$ by uniform continuity).
In relation to this we emphasize that throughout the paper what really matters is what happens
near the boundary $\po$. We also note that, while the validity
of estimate (\ref{elv}) and assumption $\DDD$ is independent of the specific distance $d(\cdot)\in\cD$ chosen, we shall
need to consider non-Euclidean distances since some of the intermediate calculations
do depend on the specific choice of distance.

We will consider operators of the form
\begin{equation}
Hu(x)=\sum_{\sscr{\mid\eta\mid = 2}
{\mid\zeta\mid = 2}}
D^{\eta}\{a_{\eta\zeta}(x)D^{\zeta}u(x)\}\; , \quad x\in\Omega\, ,
\label{eq:dr}
\end{equation}
subject to Dirichlet boundary conditions on $\po$; lower-order terms can be easily accomodated.
More precisely, we start with a matrix-valued function $a(x)=\{a_{\eta\zeta}(x)\}$ which is assumed to be
have entries in $L^{\infty}(\Omega)$ and to take its values in the set of all real,
$N(N+1)/2 \times N(N+1)/2$ matrices
($N/(N+1)/2$ is the number of multi-indices $\eta$ of length
$|\eta|= 2$). We assume that $\{a_{\eta\zeta}(x)\}$ is symmetric for all $x\in\Omega$ and
define a quadratic form $Q(\cdot)$ on the Sobolev space $H^2_0(\Omega)$ by
\[
Q(u)=\int_{\Omega}\sum_{\sscr{\mid\eta\mid=2}{\mid\zeta\mid=2}}
a_{\eta\zeta}(x)D^{\eta}u(x)D^{\zeta}u(x)\,dx,\quad u\in H^2_0(\Omega).
\]
We make the ellipticity assumption that there exist $\lambda,\Lambda>0$ such that
\[
\lambda Q_0(u) \leq Q(u)  \leq  \Lambda Q_0(u) \; , \qquad u\in H^2_0(\Omega) \, ,
\]
where $Q_0(u)=\int_{\Omega}(\Delta u)^2dx$ denotes the quadratic form corresponding to the bilaplacian
$\Delta^2$. We then define $H$ to be the associated self-adjoint operator on $L^2(\Omega)$, so that $\inprod{Hu}{u}=Q(u)$
for all $u\in\dom(H)$.

\section{Boundary decay}\label{ageneral}

Let $d(\cdot)\in\cD$. Let $\alpha >0$ be fixed and let us define
\[
\omega(x)=d(x)^{-\alpha}\; , \qquad x\in\Omega\, .
 \]
We regularize $\omega$ defining
\begin{equation}
d_n(x)=d(x)+\frac{1}{n} \quad , \qquad \omega_n(x)=d_n(x)^{-\alpha},\quad  n=1,2,\ldots\, .
\label{omega}
\end{equation}
We note that $u\in H^2_0(\Omega)$ implies $\omega_nu\in H^2_0(\Omega)$, $n\in\N$.
It is crucial for the estimates which follow that, while they contain the functions $d_n$
and $\omega_n$, they involve constants that are independent of $n\in\N$.

\begin{lemma}
Let $\alpha >0$. There exists a constant $c$ which is independent of $n\in\N$ such that
\begin{equation}
\int_{\Omega}\left( \frac{|\nabla^2u|^2}{d_n^{2\alpha}}+\frac{|\nabla u|^2}{d_n^{2+2\alpha}} +
\frac{u^2}{d_n^{4+2\alpha}}\right)dx \leq cQ(\omega_nu)\, , \qquad u\in H^2_0(\Omega).
\label{dala}
\end{equation}
\label{l2}
\end{lemma}
{\em Proof.} It suffices to prove (\ref{dala}) for all $u\in C^2_c(\Omega)$.
So let $u\in C^2_c(\Omega)$ be given and let $v=\omega_n u$, a function also in $C^2_c(\Omega)$.
Using $\DDD$ we have
\begin{eqnarray*}
\int_{\Omega}\frac{u^2}{d_n^{4+2\alpha}}dx&=&\int_{\Omega}\frac{v^2}{d_n^4}dx\\
&\leq&\int_{\Omega}\frac{v^2}{d^4}dx \\
&\leq& c\int_{\Omega}(\Delta v)^2dx \, .
\end{eqnarray*}
Similarly,
\begin{eqnarray*}
\int_{\Omega}\frac{|\nabla u|^2}{d_n^{2+2\alpha}}dx&=&
\int_{\Omega}\frac{1}{d_n^{2+2\alpha}}|\alpha d_n^{\alpha-1}v\nabla d_n +d_n^{\alpha}\nabla v|^2dx \\
&\leq& c\int_{\Omega}\frac{v^2}{d_n^4}dx +c\int_{\Omega}\frac{|\nabla v|^2}{d_n^2}dx \\
&\leq&c\int_{\Omega}(\Delta v)^2dx \, ,
\end{eqnarray*}
where we have used the fact that
\begin{equation}
\int_{\Omega}|\nabla^2v|^2dx =\int_{\Omega}(\Delta v)^2dx\, .
\label{omix}
\end{equation}
Finally, since $d$ and $d_n$ differ by a constant,
\[
u_{x_ix_j}
=d_n^{\alpha}v_{x_ix_j}+\alpha d_n^{\alpha-1}(d_{x_i}v_{x_j}+d_{x_j}v_{x_i})
+\alpha d_n^{\alpha-1}d_{x_ix_j}v +\alpha(\alpha-1)d_n^{\alpha-2}d_{x_i}d_{x_j}v,
\]
and therefore
\[
\int_{\Omega}\frac{|\nabla^2u|^2}{d_n^{2\alpha}}dx\leq 
c\biggl\{\int_{\Omega}|\nabla^2v|^2dx + \int_{\Omega}\frac{|\nabla v|^2}{d_n^2}dx +
\int_{\Omega}\frac{v^2}{d_n^4}dx  +\int_{\Omega}\frac{|\nabla^2d|^2}{d_n^2}v^2 dx \biggr\}.
\]
Since $d_n\geq d$, the second and third terms in the brackets are smaller than $c\int_{\Omega}(\Delta v)^2dx$ by 
the Hardy-Rellich inequalities $\DDD$.
The same is true for the last term by $\DD$. Thus, one more application of (\ref{omix})
concludes the proof. $\hfill //$

\begin{lemma}
Let $\alpha \in (0,1)$ and $\omega_n=d_n^{-\alpha}$. Then there exists a constant
$c>0$, independent of $n\in\N$, such that
\[
Q(u,\omega_n^2 u) \leq c \|Hu\|_2\|H^{\alpha/2}u\|_2  \, , \qquad u\in\dom(H).
\]
\label{l1}
\end{lemma}
{\em Proof.} For any $n\in\N$ and $u\in C^2_c(\Omega)$ we have
\begin{eqnarray*}
\int_{\Omega}\omega_n^{4/\alpha}u^2 dx &\leq& \int_{\Omega}\omega^{4/\alpha} u^2 dx \\
&=& \int_{\Omega}\frac{u^2}{d^4}dx \\
&\leq& cQ(u)\, .
\end{eqnarray*}
Hence $\omega_n^{4/\alpha}\leq cH$ in the quadratic form sense, which
by \cite[Lemma 4.20]{DaOPS} implies that $\omega_n^4\leq cH^{\alpha}$
(since $\alpha\in (0,1)$). Hence given $u\in\dom(H)$ we have
\begin{eqnarray*}
 Q(u,\omega_n^2 u)&\leq&\|Hu\|_2\|\omega_n^2u\|_2\\
&\leq&c\|Hu\|_2\|H^{\alpha/2}u\|_2\, ,
\end{eqnarray*}
which is the stated inequality. $\hfill //$

We can now establish a sufficient condition for the boundary decay estimates.
\begin{theorem}
Let $\alpha\in (0,1)$ be fixed and let $\omega_n=d_n^{-\alpha}$. Assume that there exist $k,k'>0$ independent of $n\in\N$ such that
\begin{equation}
Q(\omega_n u)\leq k Q(u,\omega_ n^2u) +k'\|u\|_2^2 , \qquad u\in C^2_c(\Omega),
\label{cos}
\end{equation}
for all $n\in\N$. Then there exists $c>0$ such that
\begin{equation}
\int_{\Omega}\left( \frac{|\nabla^2u|^2}{d^{2\alpha}}+\frac{|\nabla u|^2}{d^{2+2\alpha}} +
\frac{u^2}{d^{4+2\alpha}}\right)dx \leq c \|Hu\|_2\|H^{\alpha/2}u\|_2 \; , 
\quad u\in\dom(H).
\label{cosa}
\end{equation}
\label{prop:song}
\end{theorem}
{\em Proof. } The validity of (\ref{cos}) for all $u\in C^2_c(\Omega)$ implies its validity for
all $u\in H^2_0(\Omega)$ and in particular for all $u\in\dom(H)$. Hence given $u\in\dom(H)$ and applying
Lemmas \ref{l2} and \ref{l1} we conclude that there exists a constant $c$ such that for any $n\in\N$ there holds
\[
\int_{\Omega}\left( \frac{|\nabla^2u|^2}{d_n^{2\alpha}}+\frac{|\nabla u|^2}{d_n^{2+2\alpha}} +
\frac{u^2}{d_n^{4+2\alpha}}\right)dx \leq  c\biggl( \|Hu\|_2\|H^{\alpha/2}u\|_2 +\|u\|_2^2\biggr)\, .
\]
Letting $n\to +\infty$, applying the Dominated Convergence Theorem and using the fact that the spectrum
of $H$ is bounded away from zero we obtain (\ref{cosa}). $\hfill //$

%%%%%%%%%%%%%%  the property  %%%%%%%%%%%%%%%

\section{The property $\pa$}\label{theproperty}

The validity of assumption (\ref{cos}) of Theorem \ref{prop:song} will be our main interest in this section.
For the sake of simplicity, for any $\alpha\in (0,1)$ we define the property $\pa$ (relative to the distance function
$d\in\cD$) as
\[
\left\{\begin{array}{l}
 \mbox{ There exists constants $k,k'>0$ such that}\\ [0.4cm]
  \hspace{1cm}Q(d_n^{-\alpha}u)\leq k Q(u,d_n^{-2\alpha}u) +k'\|u\|_2^2 \qquad\qquad\qquad \pa \\ [0.4cm]
 \mbox{ for all $n\in\N$ and $u\in C^2_c(\Omega)$.} 
\end{array}\right.
\]
This is precisely assumption (\ref{cos}) of Theorem \ref{prop:song}.
Our aim in this section is to obtain sufficient conditions under which property $\pa$ is valid.
In the following three subsections we present three theorems that provide such conditions. The first
applies to all operators in the class under consideration; the second applies to operators of a specific type
but gives a better range of $\alpha>0$; and the third applies to small perturbations of operators in the second class.

{\bf Remark.} If $\po$ is smooth then the ground state $\phi$ of $\Delta^2$ decays as $d(x)^2$
as $x\to\po$. Hence the integral in the left-hand side of (\ref{cosa}) is not finite for $\alpha\geq 1/2$.
For this reason and throughout the rest of the paper we restrict our attention to $\alpha\in (0,1/2)$.

\subsection{General operators}

We always work in the context described at the begining of Section~\ref{ageneral}.
We recall that for $\alpha\in (0,1/2)$ we have $\omega_n=d_n^{-\alpha}=(d+\frac{1}{n})^{-\alpha}$; we also recall
that $\lambda$ and $\Lambda$ are the ellipticity constants of the operator $H$.

\begin{theorem}
There exists a computable constant $c>0$
such that property $\pa$ relative to the Euclidean distance is valid for $H$ for all $\alpha\in (0,c^{-1}\Lambda^{-1}\lambda)$.
\label{houk}
\end{theorem}
{\em Proof.} Let $u\in C^2_c(\Omega)$ be fixed. Setting $v=d_n^{-\alpha}u$ and using Leibniz' rule we have
\begin{eqnarray*}
&&\hspace{-2cm}Q(d_n^{-\alpha}u)-Q(u,d_n^{-2\alpha}u)\\
&=&Q(v)-Q(d_n^{\alpha}v,d_n^{-\alpha}v)\\
&=&\int_{\Omega}\sum_{\sscr{\mid\eta\mid=2}{\mid\zeta\mid=2}}a_{\eta\zeta}\Big( D^{\eta}vD^{\zeta}v  -
D^{\eta}(d_n^{\alpha}v)D^{\zeta}(d_n^{-\alpha}v)\Bigr)dx \\
&=&-\int_{\Omega}\sum_{\sscr{\mid\eta\mid=2}{\mid\zeta\mid=2}}\sum_{\ssscr{\gamma\leq\eta}{\delta\leq\zeta}{\gamma+\delta>0}}
c^{\eta}_{\gamma}c^{\zeta}_{\delta}a_{\eta\zeta}(D^{\gamma}d_n^{\alpha})(D^{\delta}d_n^{-\alpha})(D^{\eta-\gamma}v)(D^{\zeta-\delta}v)dx \\
&\leq&c\Lambda  \int \sum_{\sscr{0\leq i,j\leq 2}{i+j>0}}|\nabla^i d_n^{\alpha}|
\cdot |\nabla^j d_n^{-\alpha}|\cdot |\nabla^{2-i}v|\cdot |\nabla^{2-j}v| dx .
\end{eqnarray*}
But, by $\D$ and $\DD$,
\[
 |\nabla d_n^{\pm\alpha}| =\alpha d_n^{\pm\alpha -1} \; , \quad
 |\nabla^2 d_n^{\pm\alpha}|\leq c\alpha d_n^{\pm\alpha -2},
\]
and we thus obtain
\begin{eqnarray*}
Q(v)-Q(d_n^{\alpha}v,d_n^{-\alpha}v) &\leq& c\Lambda\alpha\int_{\Omega}
\biggl( |\nabla^2v|^2+\frac{|\nabla v|^2}{d_n^2} +
\frac{v^2}{d_n^4}\biggr)dx \\
&\leq& c\Lambda\lambda^{-1}\alpha Q(v).
\end{eqnarray*}
Hence, if $\alpha$ is such that $c\Lambda\lambda^{-1}\alpha<1$, then property $\pa$ is valid for $H$. $\hfill //$

%%%%%%%%%%%%%%%%%%%%%%%%%%%%%%%%
%%%%%%%%%%%%%%%%%%%%%%%%%%%%%%%%
%%%%%%%%%%%%%%%%%%%%%%%%%%%%%%%%

\subsection{Regular coefficients}\label{thebilaplacian}

The weak point of Theorem \ref{houk} is the poor information it provides on the range
of $\alpha$ for which $\pa$ is valid. In this subsection we shall consider operators of a more specific
type and for which we shall see that $\pa$ is valid for all $\alpha\in (0,1/2)$.

It will be useful in this subsection to drop the multi-index notation and write the quadratic form as
\[
Q(u)=\int_{\Omega}\sum_{i,j,k,l=1}^Na_{ijkl}u_{x_ix_j}u_{x_kx_l}dx \, , \qquad u\in H^2_0(\Omega)\, .
\]
We may clearly assume that the functions $a_{ijkl}$ have the following symmetries:
\begin{equation}
a_{ijkl}=a_{jikl}\;\; , \quad
a_{ijkl}=a_{ijlk} \;\; , \quad
a_{ijkl}=a_{klij}\, .
\label{symmetries}
\end{equation}
We make the following additional assumptions on the coefficients $\{a_{ijkl}\}$:
\begin{eqnarray}
 \ia && \mbox{There exist $\theta , \tau >0$ such that} \nonumber\\
&&\quad {\rm (a)} \qquad\mbox{each $a_{ijkl}$ is differentiable in $\{d(x)<\theta\}$}\label{kteo} \\
&&\quad {\rm (b)} \qquad |\nabla a_{ijkl}|\leq cd^{-1+\tau} \mbox{ on }\{d(x)<\theta\}  \nonumber \\
\ib &&  \sum_{i,j,k,l}a_{ijkl}(x)\xi_i\xi_k\eta_j\eta_l\leq \sum_{i,j,k,l}a_{ijkl}(x)\xi_i\xi_j\eta_k\eta_l \; , \quad \xi, \, \eta\in\R^N\, , 
\; x\in\Omega\, .
\nonumber
\end{eqnarray}
Without any loss of generality we assume that $\tau$ in $\ia$ is the same as in $\DD$.
Condition $\ib$ is a technical one, whose necessity is not clear. We present two examples in which it is valid.\nl
{\bf Example 1.} Suppose that $a_{ijkl}=b_{ij}b_{kl}$ for some non-negative $N\times N$ matrix $\{b_{ij}\}_{i,j}$.
Then $\ib$ is valid by the Cauchy-Schwarz inequality for the non-negative form $(\xi,\eta)\mapsto b_{ij}\xi_i\eta_j$.
This for example includes operators of the form $\Delta a(x)\Delta$, for which we have $a_{ijkl}=a(x)\delta_{ij}\delta_{kl}$.\nl
{\bf Example 2.} Suppose that $a_{ijkl}=\delta_{ij}\delta_{kl}a_{ik}$, where $a_{ik}=a_{ki}\geq 0$ for $i,k=1,\ldots, N$.
Then it is easily seen that $\ib$ is again valid.

We choose the distance function $d(\cdot )$ to be the one naturally associated
with $H$, that is the one induced by the Finsler metric $p(x,\eta)$ whose dual metric (cf. (\ref{dual}) below) is
\begin{equation}
p_*(x,\xi)= \Big(\sum_{i,j,k,l}a_{ijkl}(x)\xi_i\xi_j\xi_k\xi_l\Big)^{1/4}.
\label{metric}
\end{equation}
This implies in particular that the function $d(\cdot)$ satisfies
\begin{equation}
\sum_{i,j,k,l}a_{ijkl}(x)d_{x_i}d_{x_j}d_{x_k}d_{x_l}=1\, , \qquad \ale \: x\in\Omega \, .
\label{rr}
\end{equation}
Indeed, the inequality $\sum_{i,j,k,l}a_{ijkl}(x)d_{x_i}d_{x_j}d_{x_k}d_{x_l}\leq 1$ is shown in \cite[Lemma 1.3]{A}.
To prove the reverse let $y$ denote a point of differentiability of $d$. Then $y$ has a unique
nearest point $y_0\in\po$; so $d(y)=d(y,y_0)=:s$. Let $y_t$, $t\in [0,s]$, be the geodesic joining $y_0$ and $y$
parametrised by arc length so that $y_s=y$. Then for $\epsilon>0$ small we have on the one hand
\[
d(y_{s-\epsilon})-d(y)=d(y,y_{s-\epsilon})=p(y,y-y_{s-\epsilon})+o( \epsilon) ,
\]
and on the other hand, by differentiability,
\[
d(y_{s-\epsilon})-d(y)=\nabla d(y)\cdot (y_{s-\epsilon}-y) +o(\epsilon)\,.
\]
Hence
\begin{equation}
p_*(y,\nabla d(y))=\sup_{\xi\in\R^N}\frac{\nabla d(y)\cdot\xi}{p(y,\xi)}\geq \lim_{\epsilon\searrow 0}
\frac{\nabla d(y)\cdot (y_{s-\epsilon}-y)}{p(y,y-y_{s-\epsilon})}= 1 \,.
\label{dual}
\end{equation}
We note that the metric is Riemannian if the symbol of the operator $H$ is the square of a polynomial of degree two.

We assume that our basic hypotheses $\D -\DDD$ of the Introduction are valid;
Concerning in particular the validity of condition $\DD$, we note that it is satisfied if enough
regularity is imposed on the boundary and the coefficients. If for example the boundary is $C^3$
and the coefficients $a_{ijkl}$ lie in $C^3(0\leq d(x)<\theta)$, then
$d\in C^2(0\leq d(x)<\theta)$; see \cite[Section 1.3]{LN}.
On the other hand, for the Euclidean distance a $C^2$ boundary is enough \cite[p354]{GT}.

It is useful to introduce at this point a class $\cA$ of  integrals that are in a sense negligible.

{\bf Definition.} A family of quadratic integral forms $T_n(v)$, $v\in C^2_c(\Omega)$, $n\in\N$, belongs to the class $\cA$ if for any $\epsilon>0$
there exists $c_{\epsilon}>0$ (independent of $n\in\N$) such that
\begin{equation}
|T_n(v)|\leq \epsilon Q(v) +c_{\epsilon}\int_{\Omega}v^2dx\; , \qquad n\in\N, \; v\in C^2_c(\Omega).
\label{a}
\end{equation}

\begin{lemma}
Let $I_n(v)=\int b_n (D^{\gamma}v)(D^{\delta}v)dx$, $|\gamma|,|\delta|\leq 2$,
be a term that results after expanding $Q(d_n^{\alpha}v,d_n^{-\alpha}v)$ and integrating by parts a number of times. If
$b_n$ contains as a factor either a derivative of $a_{\eta\zeta}$ or a second-order derivative of $d_n$, then $(I_n)_n\in\cA$.
\label{lem:whis}
\end{lemma}
{\em Proof.} After expanding $Q(d_n^{\alpha}v,d_n^{-\alpha}v)$ (cf. (\ref{mon}) below) we obtain a linear combination of integrals,
and direct observation shows that each one of them has one of the following three forms (we switch temporarily to multi-index notation):
\begin{eqnarray*}
{\rm (a)}&& \int_{\Omega}a_{\eta\zeta}d_n^{-4+|\gamma+\delta|}(\nabla d_n)^{\eta+\zeta-\gamma-\delta}
(D^{\gamma}v)(D^{\delta}v)dx  \; , |\eta|=|\zeta|=2 \; , \; \gamma\leq\eta \; , \; \delta\leq\zeta  \, , \\
{\rm (b)}&& \int_{\Omega}a_{\eta\zeta}d_n^{-3+|\gamma|}(\nabla d_n)^{\eta-\gamma}
(D^{\zeta}d_n)v(D^{\gamma}v)dx  \; , \quad |\eta|=|\zeta|=2 \; , \; \gamma\leq\eta\, ,  \\
{\rm (c)}&& \int_{\Omega}a_{\eta\zeta}d_n^{-2}(D^{\eta}d_n)(D^{\zeta}d_n)v^2dx  \; ,\qquad  |\eta|=|\zeta|=2 \, .
\end{eqnarray*}
(These are distinguished by the number of second-order derivatives of $d_n$ that they contain -- none, one and two respectively.)
Hence all resulting integrals have the form
\[
\int_{\Omega}b_n(x)(D^{\gamma}v)(D^{\delta}v)dx\; ,  \qquad 0\leq |\gamma|,|\delta|\leq 2\, ,
\]
where $b_n(x)$ is a product of $a_{\eta\zeta}$ with powers and/or derivatives
of $d_n$ and, since $\nabla d_n$ is bounded,
\begin{equation}
|b_n(x)|\leq cd_n(x)^{-4+|\gamma+\delta|}, \qquad x\in\Omega\, .
\label{les}
\end{equation}
In cases $(b)$ and $(c)$ however, where $b_n(x)$ contains as a factor at least one second-order derivative
of $d_n$, it follows from condition $\DD$ of the Introduction that we have something more, namely
\[
|b_n(x)|\leq cd_n(x)^{-4+|\gamma+\delta|+\tau}, \qquad x\in\Omega\, .
\]
This easily implies that the integral lies in $\cA$ in this case.

Suppose now that we integrate by parts in the integral above, transfering one derivative from, say, $D^{\gamma}v$, ($|\gamma|\geq 1$), to the remaining
functions. If the derivative being transfered is $\partial/\partial x_i$, we obtain -- in an obvious notation -- the integral
\[
\int_{\Omega}\Big\{a_{\eta\zeta}d_n^{-4+|\gamma+\delta|}(\nabla d_n)^{\eta+\zeta-\gamma-\delta}(D^{\delta}v)\Big\}_{x_i}(D^{\gamma-e_i}v)dx\, .
\]
If the derivative $\partial/\partial x_i$ "hits" either $a_{\eta\zeta}$ or one of the factors that make up $(\nabla d_n)^{\!\eta+\!\zeta-\!\gamma-\!\delta}$
we obtain an integral of the form
\[ 
\int_{\Omega}b_n(x)(D^{\gamma-e_i}v)(D^{\delta}v)dx
\]
where $|b_n|\leq cd_n^{-1+\tau}d_n^{-4+|\gamma+\delta|}$; hence this integral belongs in $\cA$. $\hfill //$

{\bf Example.} We illustrate the last lemma with an example: in (\ref{pour}) below there appears the integral
\[
I_n(v)=\int_{\Omega}a_{ijkl}d_n^{-2}d_{x_i}d_{x_j}vv_{x_kx_l}dx\, .
\]
Letting $(T_n)_n,(T_n')_n$ denote elements in $\cA$ we compute
\begin{eqnarray*}
&&\int_{\Omega}a_{ijkl}d_n^{-2}d_{x_i}d_{x_j}vv_{x_kx_l}dx= \\
&=&-\int_{\Omega}(a_{ijkl})_{x_k}d_n^{-2}d_{x_i}d_{x_j}vv_{x_l}dx
+2\int_{\Omega}a_{ijkl}d_n^{-3}d_{x_i}d_{x_j}d_{x_k}vv_{x_l}dx\\
&&-\int_{\Omega}a_{ijkl}d_n^{-2}d_{x_ix_k}d_{x_j}vv_{x_l}dx
-\int_{\Omega}a_{ijkl}d_n^{-2}d_{x_i}d_{x_jx_k}vv_{x_kx_l}dx \\
&&-\int_{\Omega}a_{ijkl}d_n^{-2}d_{x_i}d_{x_j}v_{x_k}v_{x_l}dx \\
&=&\int_{\Omega}a_{ijkl}d_n^{-3}d_{x_i}d_{x_j}d_{x_k}(v^2)_{x_l}dx
-\int_{\Omega}a_{ijkl}d_n^{-2}d_{x_i}d_{x_j}v_{x_k}v_{x_l}dx +T_n(v)\\
&=&3\int_{\Omega}a_{ijkl}d_n^{-4}d_{x_i}d_{x_j}d_{x_k}d_{x_l}v^2dx
-\int_{\Omega}a_{ijkl}d_n^{-2}d_{x_i}d_{x_j}v_{x_k}v_{x_l}dx +T'_n(v)\, .\\
\end{eqnarray*}
{\bf Note.} The summation convention over repeated indices will be used from now on.
\begin{lemma}
There exists $(T_n)_n\in\cA$ such that
\begin{eqnarray}
&& Q(v)-Q(d_n^{\alpha}v,d_n^{-\alpha}v)= 2\alpha^2\int_{\Omega}a_{ijkl}d_n^{-2}d_{x_i}d_{x_j}v_{x_k}v_{x_l}dx\label{pourn}\\
&&\qquad +4\alpha^2\int_{\Omega}a_{ijkl}d_n^{-2}d_{x_i}d_{x_k}v_{x_j}v_{x_l}dx
 -(\alpha^4+11\alpha^2)\int_{\Omega}d_n^{-4}v^2dx +T_n(v)\, , 
\nonumber
\end{eqnarray}
for all $v\in C_c^2(\Omega)$.
\end{lemma}
{\em Proof.} For $\beta$ equal to $\alpha$ or to $-\alpha$ we have
\begin{equation}
(d_n^{\beta}v)_{x_ix_j}=d_n^{\beta}v_{x_ix_j} +\beta d_n^{\beta-1}d_{x_i}v_{x_j}+\beta d_n^{\beta-1}d_{x_j}v_{x_i}
+\beta(\beta-1)d_n^{\beta-2}d_{x_i}d_{x_j}v+\beta d_n^{\beta-1}d_{x_ix_j}v\, .
\label{mon}
\end{equation}
We substitute in $Q(d_n^{\alpha}v,d_n^{-\alpha}v)$ and expand. Now, by Lemma \ref{lem:whis}
all terms containing second-order derivatives of
$d_n$ belong to $\cA$. Further, the symmetries (\ref{symmetries}) of $a_{ijkl}$ give
\begin{eqnarray}
&& a_{ijkl}d_{x_i}d_{x_j}d_{x_k}v_{x_l}=a_{ijkl}d_{x_i}d_{x_k}d_{x_l}v_{x_j}=\ldots \nonumber \\
&& a_{ijkl}d_{x_i}d_{x_j}v_{x_kx_l}=a_{ijkl}d_{x_k}d_{x_l}v_{x_ix_j}=\ldots \label{ase} \\
&& a_{ijkl}d_{x_i}d_{x_l}v_{x_j}v_{x_k}=a_{ijkl}d_{x_i}d_{x_k}v_{x_j}v_{x_l}=\ldots \, .
\nonumber
\end{eqnarray}
Denoting by $(T_n)_n$ an element of $\cA$ which may change within the proof we thus arrive at
\begin{eqnarray}
Q(d_n^{\alpha}v,d_n^{-\alpha}v)&=&\int_{\Omega}a_{ijkl}\bigg\{ v_{x_ix_j}v_{x_kx_l} +2\alpha^2d_n^{-2}d_{x_i}d_{x_j}vv_{x_kx_l}
-4\alpha^2d_n^{-2}d_{x_i}d_{x_k}v_{x_j}v_{x_l}\nonumber \\
&&\qquad +4\alpha^2d_n^{-3}d_{x_i}d_{x_j}d_{x_k}vv_{x_l} +  \label{pour} \\
&&\qquad  \alpha^2(\alpha^2-1)d_n^{-4}d_{x_i}d_{x_j}d_{x_k}d_{x_l}v^2\bigg\}dx  +T_n(v).
\nonumber
\end{eqnarray}
We integrate by parts the second and fourth terms in the last integral. By Lemma \ref{lem:whis}, all terms that contain either derivatives of $a_{ijkl}$ or
second-order derivatives of $d_n$, belong in $\cA$. Hence, denoting always by $(T_n)_n$ a generic element of $\cA$ we obtain (cf. the example above)
\begin{eqnarray*}
\int_{\Omega}a_{ijkl}d_n^{-2}d_{x_i}d_{x_j}vv_{x_kx_l}dx&=&3\int_{\Omega}a_{ijkl}d_n^{-4}d_{x_i}d_{x_j}d_{x_k}d_{x_l}v^2dx\\
&&\quad -\int_{\Omega}a_{ijkl}d_n^{-2}d_{x_i}d_{x_j}v_{x_k}v_{x_l}dx +T_n(v)\,
\end{eqnarray*}
and, similarly,
\begin{eqnarray*}
\int_{\Omega}a_{ijkl} d_n^{-3}d_{x_i}d_{x_j}d_{x_k}vv_{x_l} dx &=&
\frac{3}{2}\int_{\Omega}a_{ijkl}d_n^{-4}d_{x_i}d_{x_j}d_{x_k}d_{x_l}v^2dx +T_n(v).
\end{eqnarray*}
Substituting in (\ref{pour}) yields
\begin{eqnarray*}
Q(d_n^{\alpha}v,d_n^{-\alpha}v)&=&\int_{\Omega}a_{ijkl}\bigg\{ v_{x_ix_j}v_{x_kx_l}-
4\alpha^2d_n^{-2}d_{x_i}d_{x_k}v_{x_j}v_{x_l}\\
&& -2\alpha^2d_n^{-2}d_{x_i}d_{x_j}v_{x_k}v_{x_l} 
+(\alpha^4+11\alpha^2)d_n^{-4}d_{x_i}d_{x_j}d_{x_k}d_{x_l}v^2\bigg\}dx +T_n(v) .
\end{eqnarray*}
Recalling that (\ref{rr}), relation (\ref{pourn}) follows. $\hfill //$

%%%%%%%%%%%%%%%%%%%%
\begin{lemma}
Let $v\in C^2_c(\Omega)$ and $w=d_n^{-3/2}v$. Then there holds
\begin{eqnarray*}
&&Q(v)-Q(d_n^{\alpha}v,d_n^{-\alpha}v)=
2\alpha^2\int_{\Omega}a_{ijkl}d_n d_{x_i}d_{x_j}w_{x_k}w_{x_l}dx\\
&&\qquad \qquad +4\alpha^2\int_{\Omega}a_{ijkl}d_n d_{x_i}d_{x_k}w_{x_j}w_{x_l}dx 
+(-\alpha^4+\frac{5\alpha^2}{2})\int_{\Omega}d_n^{-1}w^2dx +T_n(v),
\label{sev}
\end{eqnarray*}
where $(T_n)_n\in\cA$.
\label{lem:sun}
\end{lemma}
{\em Remark 1.} When working with the function $w$, all integrals have the form
\[
\int_{\Omega}b_n(x)(D^{\gamma}w)(D^{\delta}w)dx
\]
where the function $b_n$ satisfies
\begin{equation}
|b_n(x)|\leq cd_n(x)^{-1+|\gamma+\delta|}, \qquad x\in\Omega\, .
\label{les1}
\end{equation}
Such an integral lies in $\cA$ if in addition there holds
\[
|b_n(x)|\leq cd_n(x)^{-1+|\gamma+\delta|+\tau}, \qquad x\in\Omega\, .
\]
for some $\tau>0$; as before, these are precisely the integrals that contain
either second-order derivatives of $d_n$ or (first-order) derivatives of $a_{ijkl}$.

{\em Remark 2.} Since $d$ and $d_n$ differ by a constant, for the sake of simplicity
we shall write $d_{x_i}$ instead of $(d_n)_{x_i}$, etc.

{\em Proof of Lemma \ref{lem:sun}.} We substitute $v_{x_i}=d_n^{\frac{3}{2}}w_{x_i}+\frac{3}{2}d_n^{\frac{1}{2}}d_{x_i}w$
in (\ref{pourn}). Recalling the symmetry relations (\ref{ase}) (with $w$ in the place of $v$) and using the fact that
$a_{ijkl}d_{x_i}d_{x_j}d_{x_k}d_{x_l}=1$ we obtain
\begin{eqnarray*}
&&Q(v)-Q(d_n^{\alpha}v,d_n^{-\alpha}v)=\\
&=&4\alpha^2\int_{\Omega}a_{ijkl}d_n^{-2}d_{x_i}d_{x_k}
\Big[d_n^{\frac{3}{2}}w_{x_j}+\frac{3}{2}d_n^{\frac{1}{2}} d_{x_j}w\Big]
\Big[d_n^{\frac{3}{2}}w_{x_l}+\frac{3}{2}d_n^{\frac{1}{2}} d_{x_l}w\Big] dx +T_n(v)\\
&&\quad +2\alpha^2\int_{\Omega}a_{ijkl}d_n^{-2}d_{x_i}d_{x_j}
\Big[d_n^{\frac{3}{2}}w_{x_k}+\frac{3}{2}d_n^{\frac{1}{2}} d_{x_k}w\Big]
\Big[d_n^{\frac{3}{2}}w_{x_l}+\frac{3}{2}d_n^{\frac{1}{2}} d_{x_l}w\Big]dx\\
&&\quad -(\alpha^4+11\alpha^2)\int_{\Omega}d_n^{-1}w^2dx +T_n(v) \\
&=&4\alpha^2\int_{\Omega}a_{ijkl}d_n^{-2}d_{x_i}d_{x_k}
\Big[d_n^3w_{x_j}w_{x_l}+3d_n^2d_{x_j}w_{x_l}w+\frac{9}{4}d_nd_{x_j}d_{x_l}w^2\Big]dx\\
&&\quad +2\alpha^2\int_{\Omega}a_{ijkl}d_n^{-2}d_{x_i}d_{x_j}
\Big[d_n^3w_{x_k}w_{x_l}+3d_n^2d_{x_k}w_{x_l}w+\frac{9}{4}d_nd_{x_k}d_{x_l}w^2\Big]dx\\
&&\quad -(\alpha^4+11\alpha^2)\int_{\Omega}d_n^{-1}w^2dx +T_n(v) \\
&=&4\alpha^2\int_{\Omega}a_{ijkl}d_n d_{x_i}d_{x_k}w_{x_j}w_{x_l}dx
+2\alpha^2\int_{\Omega}a_{ijkl}d_n d_{x_i}d_{x_j}w_{x_k}w_{x_l}dx \nonumber \\
&&\quad +(-\alpha^4+\frac{5\alpha^2}{2})\int_{\Omega}d_n^{-1}w^2dx 
+18\alpha^2\int_{\Omega}a_{ijkl}d_{x_i}d_{x_j}d_{x_k}w_{x_l}wdx +T_n(v) .
\end{eqnarray*}
But the last integral belongs in $\cA$ by an integration by parts; hence the proof is complete. $\hfill //$

\begin{lemma}
Let $v\in C^2_c(\Omega)$ and $w=d_n^{-3/2}v$. Then there holds
\begin{eqnarray*}
Q(v)&=&\int_{\Omega}a_{ijkl}d_n^3w_{x_ix_j}w_{x_kx_l}dx +\frac{9}{2}\int_{\Omega}a_{ijkl}d_nd_{x_i}d_{x_j}w_{x_k}w_{x_l}dx\\
&&-3\int_{\Omega}a_{ijkl}d_nd_{x_i}d_{x_k}w_{x_j}w_{x_l}dx
+\frac{9}{16}\int_{\Omega}d_n^{-1}w^2dx +T_n(v)\, ,
\end{eqnarray*}
where $(T_n)_n$ is an element of $\cA$.
\label{lem:gold}
\end{lemma}
{\em Proof.} We have
\[
v_{x_ix_j}=d^{\frac{3}{2}}w_{x_ix_j}+\frac{3}{2}d_n^{\frac{1}{2}}d_{x_i}w_{x_j}+\frac{3}{2}d_n^{\frac{1}{2}}d_{x_j}w_{x_i}
+\frac{3}{4}d_n^{-\frac{1}{2}}d_{x_i}d_{x_j}w+\frac{3}{2}d_n^{\frac{1}{2}}d_{x_ix_j}w\, .
\]
As already mentioned, all terms involving second-order derivatives of $d_n$ belong in $\cA$.
Hence, using the symmetry relations (\ref{ase}) once more we compute,
\begin{eqnarray}
Q(v)&=&\int_{\Omega}a_{ijkl}\biggl[d_n^{\frac{3}{2}}w_{x_ix_j}+3d_n^{\frac{1}{2}}d_{x_i}w_{x_j}+\frac{3}{4}d_n^{-\frac{1}{2}}
d_{x_i}d_{x_j}w\biggr]\nonumber\\
&&\qquad\times \bigg[d_n^{\frac{3}{2}}w_{x_kx_l}+3d_n^{\frac{1}{2}}d_{x_k}w_{x_l}+\frac{3}{4}d_n^{-\frac{1}{2}}
d_{x_k}d_{x_l}w\bigg] dx +T_n(v) \label{aep}\\
&=&\int_{\Omega}a_{ijkl}\bigg[d_n^3w_{x_ix_j}w_{x_kx_l}+6d_n^2d_{x_i}w_{x_j}w_{x_kx_l}
+\frac{3}{2}d_nd_{x_i}d_{x_j}w_{x_kx_l}w\nonumber\\
&&+9d_nd_{x_i}d_{x_k}w_{x_j}w_{x_l} +\frac{9}{2}d_{x_i}d_{x_j}d_{x_k}w_{x_l}w+\frac{9}{16}d_n^{-1}d_{x_i}d_{x_j}d_{x_k}d_{x_l}w^2
\bigg]dx +T_n(v) .
\nonumber
\end{eqnarray}
The fifth term belongs in $\cA$ be a simple integration by parts. We also integrate by parts the second and third terms, obtaining respectively
\begin{eqnarray*}
\int_{\Omega}a_{ijkl}d_n^2d_{x_i}w_{x_j}w_{x_kx_l}dx&=&-2\int_{\Omega}a_{ijkl}d_nd_{x_i}d_{x_k}w_{x_j}w_{x_l}dx\\
&&\qquad +\int_{\Omega}a_{ijkl}d_nd_{x_i}d_{x_j}w_{x_k}w_{x_l}dx +T_n(v)\, ,\\
\int_{\Omega}a_{ijkl}d_nd_{x_i}d_{x_j}w_{x_kx_l}wdx &=&-\int_{\Omega}a_{ijkl}d_nd_{x_i}d_{x_j}w_{x_k}w_{x_l}dx+T_n(v)\, .
\end{eqnarray*}
Substituting in (\ref{aep}) and recalling (\ref{rr}) we obtain the stated relation. $\hfill //$

We can now prove the main theorem of this subsection. For any $\alpha\in (0,1/2)$ we define
\[
k_{\alpha}=\frac{9}{(1-4\alpha^2)(9-4\alpha^2)}.
\]
\begin{theorem}
For the operator $H$ and relative to the metric (\ref{metric})
property $\pa$ is valid for all $\alpha\in (0,1/2)$. More precisely, for any $\alpha\in (0,1/2)$
and any $k>k_{\alpha}$ there exists $k'<+\infty$ such that
\begin{equation}
Q(d_n^{-\alpha}u)\leq kQ(u,d_n^{-2\alpha}u) +k'\|u\|_2^2\, , \qquad u\in C^2_c(\Omega).
\label{mer}
\end{equation}
\label{fota}
\end{theorem}
{\em Proof.} Let $u\in C^2_c(\Omega)$ be given and let $v$ and $w$ be defined by $v=d_n^{-\alpha}u$ and $w=d_n^{-3/2}v$
respectively.
Define $\gamma_{\alpha}=(40\alpha^2-16\alpha^4)/9$ and observe that $\gamma_{\alpha}\in (0,1)$.
Applying Lemmas \ref{lem:sun} and \ref{lem:gold} and assumption (\ref{kteo}) \ib we obtain
\begin{eqnarray}
&&\hspace{-2cm}\gamma_{\alpha} Q(d_n^{-\alpha}u) -[Q(d_n^{-\alpha}u)-Q(u,d_n^{-2\alpha}u)]\nonumber\\
&=&\gamma_{\alpha} Q(v) -[Q(v)-Q(d_n^{\alpha}v,d_n^{-\alpha}v)]\nonumber\\
&=&\gamma_{\alpha}\int_{\Omega}a_{ijkl}d_n^3w_{x_ix_j}w_{x_kx_l}dx 
+\Big(\frac{9\gamma_{\alpha}}{2}-2\alpha^2\Big)\int_{\Omega}a_{ijkl}d_nd_{x_i}d_{x_j}w_{x_k}w_{x_l}dx\nonumber\\
&&\qquad -(3\gamma_{\alpha} +4\alpha^2)\int_{\Omega}a_{ijkl}d_nd_{x_i}d_{x_k}w_{x_j}w_{x_l}dx\label{con}\\
&&\qquad +\Big(\alpha^4-\frac{5\alpha^2}{2}+\frac{9\gamma_{\alpha}}{16}\Big)\int_{\Omega}a_{ijkl}d_n^{-1}d_{x_i}d_{x_j}d_{x_k}d_{x_l}w^2dx 
 +T_n(v)\nonumber\\
&\geq&\gamma_{\alpha}\int_{\Omega}a_{ijkl}d_n^3w_{x_ix_j}w_{x_kx_l}dx 
+\Big(\frac{3\gamma_{\alpha}}{2}-6\alpha^2\Big)\int_{\Omega}a_{ijkl}d_nd_{x_i}d_{x_j}w_{x_k}w_{x_l}dx\nonumber\\
&&\qquad +\Big(\alpha^4-\frac{5\alpha^2}{2}+\frac{9\gamma_{\alpha}}{16}\Big)
\int_{\Omega}a_{ijkl}d_n^{-1}d_{x_i}d_{x_j}d_{x_k}d_{x_l}w^2dx+T_n(v)\, .
\nonumber
\end{eqnarray}
Therefore
\begin{equation}
\gamma_{\alpha} Q(d_n^{-\alpha}u) -[Q(d_n^{-\alpha}u)-Q(u,d_n^{-2\alpha}u)] \geq T_n(v) \, ,
\label{fb}
\end{equation}
since the coefficient of the last integral is zero and those of the other two integrals are
non-negative.
Now, for any $\epsilon_1,\epsilon_2>0$ we have from (\ref{a})
\begin{eqnarray*}
|T_n(v)|&\leq&\epsilon_1 Q(v)+c_{\epsilon_1}\|v\|_2^2 \\
&=&\epsilon_1 Q(d_n^{-\alpha}u)+c_{\epsilon_1}\|d_n^{-\alpha}u\|_2^2 \\
&\leq&\epsilon_1 Q(d_n^{-\alpha}u)+c_{\epsilon_1}(\epsilon_2\|d_n^{-\alpha-2}u\|_2^2+c_{\epsilon_2}\|u\|_2^2) \\
&\leq&\epsilon_1 Q(d_n^{-\alpha}u)+c_{\epsilon_1}\Big(c\epsilon_2 Q(d_n^{-\alpha}u)+c_{\epsilon_2}\|u\|_2^2\Big)\, ,
\end{eqnarray*}
and therefore
\begin{equation}
|T_n(v)|\leq \epsilon Q(d_n^{-\alpha}u)+c_{\epsilon}\|u\|_2^2 
\label{con1}
\end{equation}
for any $\epsilon>0$ small. Choosing $\epsilon>0$ so that $\gamma_{\alpha}+\epsilon<1$
we obtain from (\ref{fb}) and (\ref{con1})
\[
Q(d_n^{-\alpha}u) \leq\frac{1}{1-\gamma_{\alpha}-\epsilon}Q(u,d_n^{-2\alpha}u)
 +\frac{c_{\epsilon}}{1-\gamma_{\alpha}-\epsilon}\|u\|_2^2\, .
\]
Hence (\ref{mer}) is valid with $k$ any number larger than $1/(1-\gamma_{\alpha})=k_{\alpha}$.  $\hfill //$

\subsection{Small perturbations}

In this subsection we prove a stability theorem on the validity of $\pa$. 
We denote by $M_+$ the cone of all coefficients matrices for the operators under consideration, that is
\begin{eqnarray*}
M_+&=&\biggl\{a=\{a_{\eta\zeta}\}_{|\eta|=|\zeta|=2} \; : \; a_{\eta\zeta} \mbox{ symmetric, real valued and measurable}\\ [0.2cm]
&&\qquad  \mbox{with $\lambda Q_0(u) \leq Q(u) \leq\Lambda Q_0(u)$, $u\in C^2_c(\Omega)$ $(\lambda,\Lambda>0$)}\biggr\},
\end{eqnarray*}
equipped with the uniform norm
\begin{equation}
\|a\|_{\infty}:={\rm ess\, sup}|a(x)|_{\infty}\; ;
\label{pr3}
\end{equation}
here $|a(x)|$ is the norm of the matrix $a(x)=\{a_{\eta\zeta}(x)\}_{\eta,\zeta}$ considered as an operator
on $\R^{N(N+1)/2}$.
We recall that $\lambda$, $\tilde{\lambda}$, etc, denote the lower ellipticity constant for
the operators induced by the matrices $a$, $\tilde{a}$, etc. We have

\begin{lemma}
There exists a computable constant $c>0$ such that for all $\alpha\in (0,1/2)$ there holds
\[
\int_{\Omega} |\nabla^2(d_n^{\alpha}v)|\cdot |\nabla^2(d_n^{-\alpha}v)|dx \leq cQ_0(v)\; , 
\qquad v\in C^2_c(\Omega).
\]
\label{lem:ter}
\end{lemma}
{\em Note.} For an estimate on the constant $c$ see the remark at the end of this subsection. \nl
{\em Proof.} For any $\beta\in\R$ we have
\begin{equation}
(d_n^{\beta}v)_{x_ix_j}=d^{\beta}v_{x_ix_j} +\beta d_n^{\beta-1}d_{x_i}v_{x_j}+\beta d_n^{\beta-1}d_{x_j}v_{x_i}
+\beta(\beta-1)d_n^{\beta-2}d_{x_i}d_{x_j}v+\beta d_n^{\beta-1}d_{x_ix_j}v\, .
\label{london}
\end{equation}
We write this for $\beta=\alpha$ and for $\beta=-\alpha$, and we multiply the two relations;
$d_n^{\alpha}$ cancels with $d_n^{-\alpha}\!$ and we obtain
\[
|\nabla^2(\omega_n^{-1}v)|\cdot |\nabla^2(\omega_n v)| \leq c\left\{|\nabla^2v|^2
+\frac{|\nabla v|^2}{d_n^2}+ \frac{v^2}{d_n^4}  +\frac{|\nabla^2 d_n|^2}{d_n^2}v^2\right\}.
\]
The proof is concluded by using assumption $\DD$ on $\nabla^2d$ and the Hardy-Rellich inequalities $\DDD$;
here we have also used the fact that $\int_{\Omega}|\nabla^2v|^2dx=\int_{\Omega}(\Delta v)^2dx$. $\hfill //$

\begin{prop}
Let $\alpha\in (0,1/2)$ be fixed. Assume that $\pa$ is valid for the matrix $a\in M_+$ 
relative to some distance $d(\cdot)\in\cD$ and let $k,k'>0$ be such that
\[
Q(\omega_n u) \leq k Q(u,\omega_n^2 u)+k'\|u\|_2^2\, , \quad n\in\N\, , \quad u\in C^2_c(\Omega).
\]
Then there exists a constant $c>0$
such that if $\tilde{a}\in M_+$ satisfies $\|\tilde{a}-a\|_{\infty}<\tilde{\lambda}[(1+ck)]^{-1}$, then $\pa$ is also
satisfied for $\tilde{a}$ relative to $d(\cdot)$; more precisely, there exists $\tilde{k}'<+\infty$ so that
\[
\tilde{Q}(\omega_nu)\leq \frac{k}{1-\tilde{\lambda}^{-1}(1+ck)\|a-\tilde{a}\|_{\infty}}\tilde{Q}(u,\omega_n^2 u) +
\tilde{k}'\|u\|_2^2\, , \quad n\in\N\, , \quad u\in C^2_c(\Omega).
\]
\label{prop:stability}
\end{prop}
{\em Proof.} We first note that
\begin{equation}
|\tilde{Q}(v)-Q(v)|\leq \int_{\Omega}|\tilde{a}-a|\cdot |\nabla^2v|^2dx\leq  \|\tilde{a}-a\|_{\infty}Q_0(v)\; , \qquad v\in C^2_c(\Omega).
\label{anan}
\end{equation}
Moreover, setting $v=\omega_n u$ we have from Lemma \ref{lem:ter},
\begin{eqnarray}
|\tilde{Q}(u,\omega_n^2u)-Q(u,\omega_n^2u)|&\leq& \|\tilde{a}-a\|_{\infty}\int_{\Omega}|\nabla^2u|
\cdot |\nabla^2(\omega_n^2 u)|dx\nonumber\\
&=&\|\tilde{a}-a\|_{\infty}\int_{\Omega}|\nabla^2(\omega_n^{-1}v)|\cdot |\nabla^2(\omega_n v)|dx\nonumber\\
&\leq&c \|\tilde{a}-a\|_{\infty} Q_0(v)\, .
\label{anan1}
\end{eqnarray}
From (\ref{anan}) and (\ref{anan1}) we conclude that for any $n\in\N$ and $u\in C^2_c(\Omega)$ we have
\begin{eqnarray*}
\tilde{Q}(\omega_n u)&\leq&Q(\omega_n u)+\|\tilde{a}-a\|_{\infty} Q_0(\omega_n u)\\
&\leq& kQ(u,\omega_n^2 u) +k'\|u\|_2^2 +\|\tilde{a}-a\|_{\infty} Q_0(\omega_n u)\\
&\leq&k\biggl( \tilde{Q}(u,\omega_n^2 u) +c\|\tilde{a}-a\|_{\infty}Q_0(\omega_n u)\biggr)+k'\|u\|_2^2 
 +\|\tilde{a}-a\|_{\infty} Q_0(\omega_n u)\\
&=&k\tilde{Q}(u,\omega_n^2 u)   +(1+ck)\|\tilde{a}-a\|_{\infty} Q_0(\omega_n u)+k'\|u\|_2^2\\
&\leq&k\tilde{Q}(u,\omega_n^2 u) +\tilde{\lambda}^{-1}(1+ck)\|\tilde{a}-a\|_{\infty} \tilde{Q}(\omega_n u)+k'\|u\|_2^2,
\end{eqnarray*}
from which the statement of the lemma follows. $\hfill //$

Let $\cG$ denote the cone of all coefficient matrices that satisfy assumptions $\ia$, $\ib$ of Section \ref{thebilaplacian}.
Also let $k_{\alpha}$ be as in Theorem \ref{fota}.
Combining Proposition \ref{prop:stability} and Theorem \ref{fota} we obtain immediately the following
\begin{theorem}
There exists a computable constant $c>0$ such that if for some $\alpha\in (0,1/2)$ the coefficient
matrix $a$ of the operator $H$ satisfies
\[
\dist_{L^{\infty}}(a,\cG)<\frac{\lambda}{(1+ck_{\alpha})}
\]
then $\pa$ is satisfied for $H$.
\label{grud}
\end{theorem}
{\em Proof.} Let $\tilde{a}\in\cG$ be such that
\[ 
\| a-\tilde{a}\|_{\infty}<\frac{\lambda}{1+ck_{\alpha}}
\]
By Theorem \ref{fota} $\pa$ is satisfied for $\tilde{a}$ and (\ref{mer}) is valid for any $k>k_{\alpha}$.
If in addition $k$ satisfies
\[
\| a-\tilde{a}\|_{\infty}<\frac{\lambda}{1+ck},
\]
then $\pa$ is also valid for $a$ by Proposition \ref{prop:stability}. $\hfill //$

{\bf Example.} Suppose that the coefficients $a_{\eta\zeta}$ are uniformly continuous and satisfy (\ref{kteo}) \ib. Then property $\pa$ is valid for $H$
for all $\alpha\in (0,1/2)$. This is seen by approximating $a_{\eta\zeta}$ with smooth functions using an approximate identity; note that the approximating functions
also satisfy $\ib$.

{\bf Remark.} The constant $c$ of the above proposition is precisely the constant $c$ of Lemma \ref{lem:ter}.
Precise estimates for
this constant can be easily obtained. Indeed, it follows from (\ref{london})
that for $\alpha\in (0,1/2)$ there holds modulo $\cA$
\begin{eqnarray*}
&&\int_{\Omega}|\nabla^2(d_n^{\alpha}v)|\cdot |\nabla^2(d_n^{-\alpha}v)|dx \\
&\leq& 3\int_{\Omega}\Big(
|\nabla^2v|^2 +4\alpha^2|\nabla d|^2\frac{|\nabla v|^2}{d_n^2} +\alpha^2(\alpha +1)^2|\nabla d|^4\frac{v^2}{d_n^4}\Big)dx\\
&\leq& 3\int_{\Omega}\Big(
|\nabla^2v|^2 +|\nabla d|^2\frac{|\nabla v|^2}{d_n^2} +\frac{9}{16}|\nabla d|^4\frac{v^2}{d_n^4}\Big)dx
\end{eqnarray*}
Hence, letting $c_2,c_3$ be as in $\D$, $\DDD$, we obtain (modulo $\cA$)
\[
\int_{\Omega}|\nabla^2(d_n^{\alpha}v)|\cdot |\nabla^2(d_n^{-\alpha}v)|dx 
\leq 3\Big(1+c_2^2c_3^{-1}+\frac{9}{16}c_2^4c_3^{-1}\Big)\int_{\Omega}(\Delta v)^2dx\, .
\]
In fact, since we work modulo $\cA$, the last constant can be improved
$3(1+c_2^2A^{-1}+(9/16)c_2^4B^{-1})$, where $A$ and $B$ are the {\em weak} Hardy constants, that is they satisfy
\begin{eqnarray*}
&&\int_{\Omega}|\nabla v|^2 dx \geq A
\int_{\Omega}\frac{v^2}{d^2} dx  -c'\int_{\Omega}v^2dx \\
&&\int_{\Omega}(\Delta v)^2 dx \geq B \int_{\Omega}\frac{v^2}{d^4}dx -c''\int_{\Omega}v^2dx \, .
\end{eqnarray*}
For smooth boundaries with a smooth Riemannian metric this amounts to $A=1/4$, $B=9/16$.

%%%%%%%%%%%%%%%%%%%%%%%%%%%%%%%%%%%%

\section{An application: eigenvalue stability}\label{anapplication}

In this final section we demonstrate how the boundary decay estimate of Theorem \ref{prop:song}
yield stability bounds on the eigenvalues $\{\lambda_n\}$ of $H$ under small perturbations of the boundary
$\po$. The proof follows closely the corresponding proof in \cite{Dmz} for the second-order case,
however we include it here for the sake of completeness. So we consider a distance function
$d(\cdot)\in\cD$, an operator $H$ as above
and assume that the boundary decay estimates (\ref{cosa}) are valid for some fixed $\alpha\in (0,1/2)$.
For $\epsilon>0$ we define $\Omega_{\epsilon}=\{x\in\Omega \, : \,
d(x)>\epsilon\}$. We assume that $\epsilon<\theta/2$ so that $d(x)$ is $C^2$ on
$\Omega\setminus\Omega_{2\epsilon}$.
We define $d_{\epsilon}(x)=\dist (x,\po_{\epsilon})$, $x\in\Omega$, and make
the additional assumption that there exists
$c>0$ such that for small enough $\epsilon>0$ there holds
\begin{equation}
|\nabla^2d_{\epsilon}|\leq c\mbox{ on }\{d(x)<2\epsilon\}.
\label{thola}
\end{equation}
Now, let $\tilde{\Omega}$ be a domain such that
\[ \Omega_{\epsilon}\subset \tilde{\Omega}\subset\Omega\, ;\]
we do not make any regularity assumptions on $\partial\tilde{\Omega}$.
We denote by $\{\tilde{\lambda}_n\}$ the eigenvalues of the operator
$\tilde{H}$ on $L^2(\tilde{\Omega})$, which is defined
by restricting the quadratic form $Q(\cdot)$ on $H^2_0(\tilde{\Omega})$.

Let $\phi$ be a
non-negative, smooth, increasing function on $\R$ such that
\[ \phi(t)=\darr{0,}{t\leq 0,}{1,}{t\geq 1.}\]
We define a $C^2$ cut-off function $\tau$ on $\Omega$ by
\[
\tau(x)=\darr{0,}{x\in\Omega\setminus\Omega_{\epsilon}\, ,}{\phi({d_{\epsilon}}(x)/\epsilon),}{x\in\Omega_{\epsilon}\, .}
\]
Note that $\tau(x)=1$ when $d(x)>2\epsilon$; moreover (\ref{thola}) yields
\[
|\tau (x)|\leq 1 \quad , \quad |\nabla\tau(x)|\leq c\epsilon^{-1} \quad , \quad |\nabla^2\tau(x)|\leq c\epsilon^{-2}\, .
\]

%%%%%%%%%%%%%%%%%%%%%%%%%%%%%%%%%%%%%%%%%%%%%%%%%%%

Let us now denote by $\{\phi_n\}$ the normalized eigenfunctions of $H$.
For $n\geq 1$ we set
\[
L_n={\rm lin}\{\phi_1,\ldots ,\phi_n\} \; , \quad  \tilde{L}_n={\rm lin}\{\tau\phi_1,\ldots ,\tau \phi_n\},
\]
and observe that $\tilde{L}_n\subset H^2_0(\tilde{\Omega})$.
\begin{lemma}
There exists a constant $c>0$ such that for all $\epsilon>0$ small
and all $u\in\dom(H)$ there holds
\begin{eqnarray*}
\ia&& |Q(\tau u)-Q(u)|\leq
c\epsilon^{2\alpha} \|Hu\|_2\|H^{\alpha/2}u\|_2  \, , \\ [0.2cm]
\ib && \Bigl| \|\tau u\|_2 -\|u\|_2 \Bigr| \leq c\epsilon^{2+\alpha}   \|Hu\|_2^{1/2}\|H^{\alpha/2}u\|_2^{1/2} \, .
\end{eqnarray*}
\label{oli}
\end{lemma}
{\em Proof.} Let $u\in\dom(H)$. On $\Omega_{2\epsilon}$ we have $\tau(x)=1$, hence
\begin{eqnarray*}
|Q(\tau u)-Q(u)| &=&\bigg|\int_{\Omega}
\sum_{\sscr{\mid\eta\mid=2}{\mid\zeta\mid=2}}a_{\eta\zeta}
\left\{  (D^{\eta}(\tau u))(D^{\zeta}(\tau u)) - (D^{\eta}u)(D^{\zeta}u)
\right\}dx\bigg| \\
&\leq&c\int_{d(x)<2\epsilon}
 \Bigl\{ |\nabla^2(\tau u)|^2  +|\nabla^2 u|^2\Bigr\} dx\\
&\leq&c\int_{d(x)<2\epsilon}
 \Bigl\{  |\nabla^2 u|^2 +|\nabla\tau|^2|\nabla u|^2  +|\nabla^2\tau|^2u^2\Bigr\} dx\\
&\leq&c\int_{d(x)<2\epsilon}
 \Bigl\{  |\nabla^2 u|^2 +\epsilon^{-2} |\nabla u|^2 +\epsilon^{-4} u^2 \Bigr\} dx\\
&\leq&c\epsilon^{2\alpha}\int_{d(x)<2\epsilon}
 \Bigl\{ \frac{|\nabla^2 u|^2}{d^{2\alpha}} +\frac{|\nabla u|^2}{d^{2+2\alpha}} +\frac{|u|^2}{d^{4+2\alpha}} \Bigr\}dx,
\end{eqnarray*}
from which $\ia$ follows by means of Theorem \ref{prop:song}. Similarly,
\begin{eqnarray*}
\Bigl| \|\tau u\|_2 -\|u\|_2 \Bigr|^2 &\leq& \|\tau u -u\|_2^2 \\
&\leq&\int_{d(x)<2\epsilon} |u|^2dx \\
&\leq&\epsilon^{4+2\alpha}\int_{d(x)<2\epsilon} \frac{u^2}{d^{4+2\alpha}}dx\\
&\leq&c\epsilon^{4+2\alpha}  \|Hu\|_2\|H^{\alpha/2}u\|_2 \, ,
\end{eqnarray*}
from which $\ib$ follows. $\hfill //$

%%%%%%%%%%%%%%%%%%%%%%%%%%%%%%%%%%%%%%%%%

\begin{theorem}
Assume that there exists a distance
function $d\in\cD$ and an $\alpha\in (0,1/2)$ such that (\ref{cosa}) is satisfied. Assume also
that (\ref{thola}) is valid.
Then there exists $c,c'>0$ such that for each $n\geq 1$ there holds
\begin{equation}
0<\lambda_n \leq\tilde{\lambda}_n\leq\lambda_n +c\lambda_n^{5/4}\epsilon^{2\alpha}\, ,
\label{diag}
\end{equation}
for all $\epsilon>0$ satisfying $\epsilon^{2\alpha}<c'\lambda_n^{-5/4}$.
\label{leve}
\end{theorem}
{\em Proof.} We fix $n\geq 1$. Since $\tilde{L}_n\subset H^2_0(\tilde{\Omega})$ we have by min-max
\begin{eqnarray}
\tilde{\lambda}_n&\leq&\sup\{ Q(v)/\|v\|_2^2 \; : \; v\in \tilde{L}_n\} \nonumber\\
&=&\sup\{ Q(\tau u)/\|\tau u\|_2^2 \; : \; u\in L_n\}.
\label{param}
\end{eqnarray}
Now, let $u\in L_n$ be given. It follows from Lemma \ref{oli} $\ia$ that
\begin{eqnarray}
Q(\tau u)&\leq& Q(u)+c\epsilon^{2\alpha}\|Hu\|_2\|H^{\alpha/2}u\|_2  \nonumber\\
&\leq&Q(u)+c\epsilon^{2\alpha}\lambda_n^{5/4}\|u\|_2^2 \, . \label{nice}
\end{eqnarray}
Similarly Lemma \ref{oli} $\ib$ gives
\begin{eqnarray}
\|\tau u\|_2^2 &\geq& \|u\|_2^2 - c\epsilon^{2+\alpha} \lambda_n^{5/8}\|u\|_2
(\|u\|_2+\|\tau u\|_2)\nonumber \\
&\geq& \|u\|_2^2 -c\epsilon^{2+\alpha} \lambda_n^{5/8}\|u\|_2^2 \, .
\label{nice1}
\end{eqnarray}
Assumimg in addition that $\|u\|_2=1$ we thus obtain from (\ref{nice}) and (\ref{nice1})
that
\begin{eqnarray*}
\frac{Q(\tau u)}{\|\tau u\|_2^2}&\leq&\frac{Q(u)+
c\epsilon^{2\alpha}\lambda_n^{5/4}}{1 -c\epsilon^{2+\alpha}\lambda_n^{5/8}}\\
&\leq& Q(u) +c\lambda_n^{5/4}\epsilon^{2\alpha}\\
&\leq&\lambda_n +c\lambda_n^{5/4}\epsilon^{2\alpha},
\end{eqnarray*}
where for the second inequality we have used the fact that $\epsilon^{2\alpha}<c'\lambda_n^{-5/4}$,
with $c'$ small enough but fixed (and independent of $n$ and $\epsilon$).
Hence (\ref{param}) implies
\[
\lambda_n \leq \tilde{\lambda}_n\leq\lambda_n +c\lambda_n^{5/4}\epsilon^{2\alpha},
\]
which completes the proof of the theorem. $\hfill //$

{\bf Remark.} In the case where $\Omega=B(1)$ and $\tilde{\Omega}=B(1-\epsilon)$ we have $\tilde{\lambda}_n
=(1-\epsilon)^{-4}\lambda_n$ and hence $\tilde{\lambda}_n-\lambda_n=4\lambda_n(\epsilon+O(\epsilon^2))$.
Hence the value $\alpha=1/2$ is the best possible for estimate (\ref{diag}).

{\bf Acknowledgment} I thank Professor Z. Shen for bringing to my attention the article \cite{LN}.
Partial  support by the RTN European network Fronts--Singularities HPRN-CT-2002-00274 is also acknowledged.

%%%%%%%%%%%%%%%%%%%%%

%&&&&&&&&&&&&&&&&&&&&&&&&&&&&&&&&&&&&&&&&&&&&&&&&&&&&&&&&

%\newpage


\begin{thebibliography}{RRR}

\bibitem[A]{A}{Agmon S. Lectures on exponential decay of solutions of second-order elliptic equations.
Mathematical Notes, Princeton University Press, 1984.}

\bibitem[B]{B}{Barbatis G. Explicit estimates on the fundamental solution of higher-order
parabolic equations with measurable coefficients. {\em J. Differential Equations}{ \bf 174} (2001) 442-463.}

\bibitem[BFT]{BFT}{Barbatis G., Filippas S. and Tertikas A. Refined geometric $L^p$ Hardy
inequalities. {\em Commun. Contemp. Math.} {\bf 5} (2003) 869-881.}

\bibitem[BT]{BT}{Barbatis G. and Tertikas A. On a class of Rellich inequalities, {\em J. Comp. Appl. Math.}
{\bf 194} (2006) 156-172.}

\bibitem[D1]{DaOPS}{Davies E.B. One parameter semigroups. Academic
Press, 1980.}

\bibitem[D2]{Dmz}{Davies E.B. Eigenvalue stability bounds via weighted
Sobolev spaces. {\em Math. Z.} {\bf 214} (1993) 357-371.}

\bibitem[D3]{Dqjo}{Davies E.B. The Hardy constant. {\em Quart. J. Math. Oxford} (2) {\bf 46}
(1995) 417-431.}

\bibitem[D4]{Dcam}{Davies E.B. Sharp boundary estimates for elliptic operators. {\em
Math. Proc. Cambridge Philos. Soc.} {\bf 129} (2000) 165-178.}

\bibitem[EHK]{EHK}{Evans W.D., Harris D.J. and Kauffman R.M. Boundary behaviour of Dirichlet
eigenfunctions of second order elliptic operators. {\em Math. Z.} {\bf 204} (1990) 85-115.}

\bibitem[FHT]{FHT}{Fleckinger J., Harrell E.M. II and de Th\'{e}lin F. Boundary behavior
and estimates for solutions of equations containing the $p$-Laplacian. {\em Electron. J.
Differential Equations} {\bf 38} (1999) 19pp.}

\bibitem[GT]{GT}{Gilbarg D. and Trudinger N.S. Elliptic partial differential equations of
second order. Springer, 1983.}

\bibitem[LN]{LN}{Li Y.Y. and Nirenberg L. The distance function to the boundary, Finsler geometry
and the singular set of viscosity solutions of some Hamilton-Jacobi equations.
{\em Comm. Pure Appl. Math.} {\bf 58} (2005) 85-146.}

\bibitem[M]{M}{Mason C. Perturbation of domain: singular Riemannian metrics.
{\em Proc. London Math. Soc.} {\bf 84} (2002) 473-491.}

\end{thebibliography}
\end{document}